\documentclass[11pt]{article}
\usepackage{amsmath}
\usepackage{amssymb}
\usepackage{latexsym}
\begin{document}
\title{\sc \bf {A note on relations between Hom-Malcev algebras and Hom-Lie-Yamaguti algebras}}
\author{{\bf Donatien Gaparayi}{\footnote{E-mail address: gapadona@yahoo.fr}} \\ 
Ecole Normale Sup\'erieure, BP 6983 Bujumbura, Burundi.\\ and \\ {\bf A. Nourou Issa}{\footnote{E-mail address: woraniss@yahoo.fr}} \\ D\'epartement de Math\'ematiques, 
Universit\'e d'Abomey-Calavi, \\ 01 BP 4521 Cotonou 01, B\'enin}
\date{}
\maketitle
\begin{abstract}
A Hom-Lie-Yamaguti algebra, whose ternary operation expresses through its binary one in a specific way, is a multiplicative Hom-Malcev 
algebra. Any multiplicative Hom-Malcev algebra over a field of characteristic zero has a natural Hom-Lie-Yamaguti structure. 
\end{abstract}
\vspace{1cm}
{\bf 1 Introduction and results} \\
\par
Motivated by quasi-deformations of Lie algebras of vector fields, including q-deformations of Witt and Virasoro algebras, Hom-Lie algebras were introduced in 
\cite{HLS}. The intensive development of the theory of Hom-algebras started by the introduction of Hom-associative algebras in \cite{MS1}, where it is shown that the 
commutator algebra of a given Hom-associative algebra is a Hom-Lie algebra. Since then, various types of Hom-algebras were introduced and investigated (see, 
e.g., \cite{AMS}, \cite{AI}, \cite{GI1}, \cite{Issa1}, \cite{Makh}-\cite{MS2}, \cite{Yau3}-\cite{Yau5} and references therein). We refer to \cite{HLS}, \cite{Makh}, 
\cite{MS1}, \cite{Yau1}, \cite{Yau2} for fundamentals on Hom-algebras.
\par
A rough description of a Hom-type algebra is that it is a generalization of a given type of algebras by twisting its defining identity (identities) by a 
linear self-map in such a way that, when the twisting map is the identity map, one recovers the original type of algebras (in \cite{Yau2} a general strategy to 
twist a given algebraic structure into its corresponding Hom-type via an endomorphism is given). Following this pattern, a number of binary, $n$-ary or 
binary-ternary algebras are twisted into their Hom-version. For instance, by twisting alternative algebras, Hom-alternative algebras are introduced in \cite{Makh}, 
and the notion of a Hom-Malcev algebra is introduced in \cite{Yau3}, where it is shown that the commutator Hom-algebra of a Hom-alternative algebra is a Hom-Malcev 
algebra (this is the Hom-analogue of the Malcev's result stating that alternative algebras are Malcev-admissible). In fact, Malcev algebras were introduced 
by A.I. Mal'tsev in \cite{Mal} (calling them Moufang-Lie algebras) as commutator algebras of 
alternative algebras and also as tangent algebras to local smooth Moufang loops. The terminology ``Malcev algebra'' is introduced in \cite{Sag}, where a 
systematic study of such algebras is undertaken.
\par
All vector spaces and algebras throughout will be over a ground field $\mathbb K$ of characteristic $0$.
\par
A {\it Malcev algebra} \cite{Mal} is a nonassociative algebra $A$ with an anticommutative binary operation ``$[,]$'' satisfying the {\it Malcev identity} \\
\par
$J(x,y,[x,z]) = [J(x,y,z), x]$ \hfill (1) \\
\\
for all $x,y,z$ in $A$, where $J(x,y,z) = {\circlearrowleft}_{(x,y,z)} [[x,y],z]$ denotes the {\it Jacobian} and $ {\circlearrowleft}_{(x,y,z)} $ means the 
sum over cyclic permutation of $x,y,z$. A {\it Hom-Malcev algebra} \cite{Yau3} is a {\it Hom-algebra} $(A, [,], \alpha )$ such that its binary operation ``$[,]$`` 
is anticommutative and that the {\it Hom-Malcev identity} \\
\par
$J_{\alpha}(\alpha (x), \alpha (y), [x,z]) = [J_{\alpha}(x,y,z), {\alpha}^{2} (x)]$ \hfill (2) \\
\\
holds for all $x,y,z$ in $A$, where $J_{\alpha}(x,y,z) = {\circlearrowleft}_{(x,y,z)} [[x,y], \alpha (z)]$ is the {\it Hom-Jacobian}. Observe that 
$J_{\alpha}(x,y,z)$ is completely skew-symmetric in its three variables and when the {\it twisting map} $\alpha$ is the identity map, $\alpha = id$, then 
the Hom-Malcev algebra $(A, [,], \alpha )$ reduces to the Malcev algebra $(A, [,])$. The Hom-Malcev algebra $(A, [,], \alpha )$ is said to be 
{\it multiplicative} if $\alpha ([x,y]) = [\alpha (x), \alpha (y)]$, for all $x,y$ in $A$. We assume in this paper that all Hom-algebras are multiplicative.
In \cite{Issa2} an identity, equivalent to the Hom-Malcev identity (2), is pointed out (see also section 2).
\par
As for binary algebras, $n$-ary algebras can be twisted into their Hom-version (see, e.g., \cite{AMS}, \cite{Yau4}, \cite{Yau5}). Other interesting types of algebras 
are binary-ternary algebras, i.e. algebras with one (or more) binary operation and one (or more) ternary operation. A well-known class of such algebras is the
one of Lie-Yamaguti algebras which are introduced by K. Yamaguti \cite{Yam1} as a generalization of Lie triple systems (this motivates the name 
''generalized Lie triple systems`` used in \cite{Yam1} for these algebras; in \cite{Kik} they are called ''Lie triple algebras`` and the terminology ''Lie-Yamaguti 
algebras`` is introduced in \cite{KW} to call these algebras). It turns out that the operations of a Lie-Yamaguti algebra characterize the torsion and 
curvature tensors of the Nomizu's canonical connection on reductive homogeneous spaces \cite{Nom}.
\par
A {\it Lie-Yamaguti algebra} ({\it LYA}) $(A,*, \{ , , \})$ is a vector space $A$ together with a binary operation $*: A \times A \rightarrow A$ and a 
ternary operation $\{ , , \} : A \times A \times A \rightarrow A$ such that \\
\par
{\bf (LY1)} $x*y = - y*x$,
\par
{\bf (LY2)} $\{x,y,z\} = - \{y,x,z\}$,
\par
{\bf (LY3)} ${\circlearrowleft}_{(x,y,z)} [(x*y) * z + \{x,y,z\}] = 0$,
\par
{\bf (LY4)} ${\circlearrowleft}_{(x,y,z)} \{ x*y , z, u \} = 0$,
\par
{\bf (LY5)} $\{x, y, u*v \} = \{x,y,u\}*v + u*\{x,y,v\}$,
\par
{\bf (LY6)} $ \{x,y, \{u,v,w \} \} = \{ \{ x,y,u \},v,w \} + \{u, \{ x,y,v \},w \}$
\par
\hspace{4cm} $+ \{ u,v, \{x,y,w\} \}$, \\
for all $u,v,w,x,y,z$ in $A$. \\
\par
One observes that if $x*y = 0$, for all $x,y$ in $A$, then $(A,*, \{ , , \})$ reduces to a {\it Lie triple system} $(A,\{ , , \})$ and if $\{x,y,z\} = 0$ for 
all $x,y,z$ in $A$ then one gets a Lie algebra $(A,*)$. Motivated by recent developments of the theory of Hom-algebras, Theorem 2.3 in \cite{Yau2} is extended to 
the study of a twisted deformation of {\it Akivis algebras} (see references in \cite{Issa1}) which constitutes a very general class of binary-ternary algebras. It is 
proved (\cite{Issa1}, Corollary 4.5) that every Akivis algebra $A$ can be twisted into a {\it Hom-Akivis algebra} via an endomorphism of $A$ (this is the first
extension of Theorem 2.3 in \cite{Yau2} to the category of binary-ternary Hom-algebras). Following this line, Hom-Lie-Yamaguti algebras are introduced in 
\cite{GI1} as a twisted generalization of Lie-Yamaguti algebras (the cohomology theory and representation theory of Hom-Lie-Yamaguti algebras are recently
developed in \cite{Ma} and \cite{Zh}).
\par
A {\it Hom-Lie-Yamaguti algebra} (Hom-LYA for short) \cite{GI1} is a quadruple \\ $(A,*, \{ , , \}, \alpha)$ in which $A$ is a $\mathbb K$-vector space, 
``$*$'' a binary operation and ``$\{ , , \}$'' a ternary operation on $A$, and $\alpha : A \rightarrow A$ a linear map such that \\
\par
{\bf (HLY1)} $\alpha ( x * y) = \alpha (x) * \alpha (y)$,
\par
{\bf (HLY2)} $\alpha (\{x,y,z\}) = \{\alpha (x), \alpha (y), \alpha (z) \}$,
\par
{\bf (HLY3)} $x*y = - y*x$,
\par
{\bf (HLY4)} $\{x,y,z\} = - \{y,x,z\}$,
\par
{\bf (HLY5)} ${\circlearrowleft}_{(x,y,z)} [(x*y) * \alpha (z) + \{x,y,z\}] = 0$,
\par
{\bf (HLY6)} ${\circlearrowleft}_{(x,y,z)} \{ x*y , \alpha (z), \alpha (u) \} = 0$,
\par
{\bf (HLY7)} $\{\alpha (x), \alpha (y), u*v \} = \{x,y,u\}* {\alpha}^{2}(v) + {\alpha}^{2}(u)*\{x,y,v\}$,
\par
{\bf (HLY8)} $ \{ {\alpha}^{2}(x), {\alpha}^{2}(y), \{ u,v, w \} \} = \{ \{ x,y,u \}, {\alpha}^{2}(v), {\alpha}^{2}(w) \}$
\par
\hspace{6.0cm} $+ \{ {\alpha}^{2}(u), \{ x,y,v \}, {\alpha}^{2}(w) \}$
\par
\hspace{6.0cm} $+ \{ {\alpha}^{2}(u), {\alpha}^{2}(v), \{x,y,w\} \}$, \\
for all $u,v,w,x,y,z$ in $A$.
\par
Note that the conditions {\bf (HLY1)} and {\bf (HLY2)} mean the multiplicativity of $(A,*, \{ , , \}, \alpha)$. Examples of Hom-LYA could be found in
\cite{GI1}, \cite{GI2}.\\
\par
{\bf Remark.} (i) If $\alpha = Id$, then the Hom-LYA $(A,*, \{ , , \}, \alpha)$ reduces to a LYA $(A,*, \{ , , \})$ (see {\bf (LY1)}-{\bf (LY6)}).
\par
(ii) If $x*y = 0$, for all $x,y \in A$, then $(A,*, \{ , , \}, \alpha)$ is a multiplicative Hom-Lie triple system $(A,\{ , , \}, {\alpha}^{2})$ 
and, subsequently, a multiplicative ternary Hom-Nambu algebra since any Hom-Lie triple system is automatically a ternary Hom-Nambu algebra (see 
\cite{Yau4} for Hom-Lie triple systems and \cite{AMS} for Hom-Nambu algebras).
\par
(iii) If $\{x,y,z\} = 0$ for all $x,y,z \in A$, then the Hom-LYA $(A,*, \{ , , \}, \alpha)$ becomes a Hom-Lie algebra $(A,*, \alpha)$. \\
\par
It is shown (\cite{GI1}, Corollary 3.2) that every LYA $(A,*, \{ , , \})$ can be twisted into a Hom-LYA via an endomorphism of $(A,*, \{ , , \})$.
\par
The relationships between LYA and Malcev algebras are investigated by K. Yamaguti in \cite{Yam2}, \cite{Yam3}. In \cite{Yam2} (Theorem 1.1), relaying on 
a result in \cite{Sag} (Proposition 8.3), it is proved that the Malcev identity is equivalent to {\bf (LY5)} in an anticommutative algebra over a field of 
characteristic not 2 or 3 with ``$\{ , , \} $'' defined in a specific way. Moreover, 
any Malcev algebra over a field of characteristic not 2 has a natural LYA structure (\cite{Yam2}, proof of Theorem 2.1). Besides, when the ternary operation 
of a given LYA expresses in a specific way through its binary one, then such a LYA reduces to a Malcev algebra (\cite{Yam3}, Theorem 1.1).
\par
The purpose of this note is the study of the twisted version of K. Yamaguti's results relating Malcev algebras and LYA (\cite{Yam2}, \cite{Yam3}) that is, in a 
similar way, we shall relate Hom-Malcev algebras and Hom-LYA. We stress that although the analogue of Yamaguti's results are shown below to hold in the 
Hom-algebra setting, the methods used in the proofs of these results still cannot be reported in the case of Hom-algebras. So we proceed otherwise as it could
be seen in what follows.
\par 
Our investigations are based on the trilinear composition \\
\par
$ \{x,y,z\} := xy* \alpha (z) - yz* \alpha (x) - zx* \alpha (y) $, \hfill (3) \\
\\
where ``$*$'' will denote the binary operation of either the given Hom-Malcev algebra or the Hom-LYA and juxtaposition is used in order to reduce the 
number of braces i.e., e.g., $xy* \alpha (z)$ means $(x*y)* \alpha (z)$. We shall prove: \\
\par
{\bf Theorem 1.1.} {\it Let $(A,*, \{ , , \}, \alpha)$ be a Hom-LYA. If its ternary operation ``$\{ , , \}$'' expresses through its binary one ``$*$'' 
as in} (3) {\it for all $x,y,z$ in $A$, then $(A,*, \alpha)$ is a multiplicative Hom-Malcev algebra}. \\
\par
{\bf Theorem 1.2.} {\it Let $(A,*, \alpha)$ be a multiplicative Hom-Malcev algebra. If define on $(A,*, \alpha)$ a ternary operation by} (3), {\it then $A$ 
has a Hom-LYA structure}. \\
\par
The next section is devoted to the proofs of Theorems 1.1 and 1.2. Some other results are also mentioned. \\
\par
{\bf 2 Proofs} \\
\par
In \cite{Issa2} it is proved that, in an anticommutative Hom-algebra $(A, [,] , \alpha)$, the Hom-Malcev identity (2) is equivalent to the identity \\
\par
$J_{\alpha}(\alpha (x), \alpha (y) ,[u,v]) = [J_{\alpha}(x,y,u),{\alpha}^{2}(v)] + [{\alpha}^{2}(u), J_{\alpha}(x,y,v)]$
\par
\hspace{3.5cm} $ -2 J_{\alpha}(\alpha (u), \alpha (v) ,[x,y])$ \hfill (4) \\
\\
so that (4) can also be taken as a defining identity of Hom-Malcev algebras. Now we write (3) in an equivalent suitable form as \\
\par
$ \{ x,y,z \} = -J_{\alpha}(x,y,z) + 2 xy * \alpha (z)$. \hfill (5) \\
\\
{\bf Proof of Theorem 1.1.} Observe that (5) and multiplicativity imply \\
\par
$\{ \alpha (x), \alpha (y), z \} = -J_{\alpha}(\alpha (x), \alpha (y),z) + 2 \alpha (xy*z)$. \hfill (6) \\
\\
Then, putting (6) in {\bf (HLY7)}, we get
\par
$-J_{\alpha}(\alpha (x), \alpha (y) , u*v) = - J_{\alpha}(x,y,u)*{\alpha}^{2}(v) + {\alpha}^{2}(u)* (-J_{\alpha}(x,y,v)) $
\par
\hspace{4cm} $ + (2xy * \alpha (u)) * {\alpha}^{2}(v) + {\alpha}^{2}(u) * (2xy* \alpha (v))$
\par
\hspace{4cm}
$ - 2 \alpha (xy*uv)$ \\
and this last equality is written as \\
\par
$J_{\alpha}(\alpha (x), \alpha (y) ,u*v) = J_{\alpha}(x,y,u)*{\alpha}^{2}(v) + {\alpha}^{2}(u)* J_{\alpha}(x,y,v)$
\par
\hspace{3.5cm} $ -2 J_{\alpha}(\alpha (u), \alpha (v) ,x*y)$, \\
\\
which is (4). Therefore $(A, * , \alpha)$ is a Hom-Malcev algebra. \hfill $\square$ \\
\par
{\bf Remark.} For $\alpha = Id$, the ternary operation (5) reduces to the ternary operation, defined by the relation (1.4) in \cite{Yam2}, 
that is considered in Malcev algebras. Thus Theorem 1.1 above is the Hom-analogue of the result of K. Yamaguti \cite{Yam2}, which is the converse of a result 
of A.A. Sagle (\cite{Sag}, Proposition 8.3). The Hom-version of the Sagle's result is the following \\
\par
{\bf Proposition 2.1} {\it Let $(A, * , \alpha)$ be a multiplicative Hom-Malcev algebra and define on $(A, * , \alpha)$ a ternary operation by} (5). {\it Then} \\
\par
$\{ \alpha (x), \alpha (y),u*v \} = \{ x,y,u \}* {\alpha}^{2}(v) + {\alpha}^{2}(u)*\{ x,y,v \}$ \hfill (7) \\
\\
{\it for all $u,v,x,y$ in $A$}. \\
\par
{\bf Proof.} We write the identity (4) as
\par
$-J_{\alpha}(\alpha (x), \alpha (y) , u*v) = - J_{\alpha}(x,y,u)*{\alpha}^{2}(v) + {\alpha}^{2}(u)* (-J_{\alpha}(x,y,v)) $
\par
\hspace{4.5cm}$+ 2 J_{\alpha}(\alpha (u), \alpha (v) ,x*y)$ \\
i.e.
\par
$-J_{\alpha}(\alpha (x), \alpha (y) ,u*v) = - J_{\alpha}(x,y,u)*{\alpha}^{2}(v) + {\alpha}^{2}(u)* (-J_{\alpha}(x,y,v)) $
\par
\hspace{4.5cm}$+ 2 \alpha (u*v) * \alpha (x*y) + 2(\alpha (v) * xy) * {\alpha}^{2}(u)$
\par
\hspace{4.5cm}$+ 2(xy * \alpha (u)) * {\alpha}^{2}(v)$ \\
or
\par
$-J_{\alpha}(\alpha (x), \alpha (y) , u*v) + 2 \alpha (x*y) * \alpha (u*v)$
\par
\hspace{4cm}$ = (-J_{\alpha}(x,y,u) + 2(xy* \alpha (u))) * {\alpha}^{2}(v) $
\par
\hspace{4.5cm}$+ {\alpha}^{2}(u) * (-J_{\alpha}(x,y,v) + 2(xy * \alpha (v)))$. \\
This last equality (according to (5) and using multiplicativity) means that
\par
$\{ \alpha (x), \alpha (y), u*v \} = \{ x,y,u \} * {\alpha}^{2}(v) + {\alpha}^{2}(u) * \{ x,y,v \} $ \\
and therefore the proposition is proved. \hfill $\square$ \\
\par
Observe that (7) is just {\bf (HLY7)} in the definition of a Hom-LY algebra. Combining Theorem 1.1 and Proposition 2.1, we get the following \\
\par
{\bf Corollary 2.2.} {\it In an anticommutative Hom-algebra $(A, * , \alpha)$, the Hom-Malcev identity} (2) {\it is equivalent to} (7), {\it with 
``$\{ , , \}$'' defined by} (5). \hfill $\square$ \\
\par
The untwisted counterpart of Corollary 2.2 is Theorem 1.1 in \cite{Yam2}. \\
\par
{\bf Proof of Theorem 1.2.} We must prove the validity in $(A, * , \alpha)$ of the set of identities {\bf (HLY1)}-{\bf (HLY8)}. In the transformations below,
we shall use the complete skew-symmetry of the Hom-Jacobian $J_{\alpha}(x,y,z)$ in $(A, * , \alpha)$.
\par
The multiplicativity of $(A, * , \alpha)$ implies {\bf (HLY1)} and {\bf (HLY2)}. The skew-symmetry of ``$*$''  is {\bf (HLY3)} and it implies 
$ \{ x,y,z \} = - \{ y,x,z \}$ which is {\bf (HLY4)}. Next,
\begin{eqnarray*}
J_{\alpha}(x,y,z) &+& {\circlearrowleft}_{(x,y,z)}\{ x,y,z \} = J_{\alpha}(x,y,z) - J_{\alpha}(x,y,z) + 2xy*\alpha (z) \nonumber \\
&-& J_{\alpha}(y,z,x) + 2yz*\alpha (x) - J_{\alpha}(z,x,y) + 2zx*\alpha (y) \nonumber \\
&=& -J_{\alpha}(y,z,x) - J_{\alpha}(z,x,y) +  2xy*\alpha (z) + 2yz*\alpha (x) \nonumber \\
&+& 2zx*\alpha (y) \nonumber \\
&=& -2J_{\alpha}(x,y,z) + 2J_{\alpha}(x,y,z) \nonumber \\
&=& 0 \nonumber
\end{eqnarray*}
so we get {\bf (HLY5)}. Now  consider $ {\circlearrowleft}_{(x,y,z)} \{ x*y, \alpha (z),  \alpha (u)\}$ and note that, by (4), we have
\begin{eqnarray*}
& & \{ x*y, \alpha (z),  \alpha (u)\} = - J_{\alpha}(x*y, \alpha (z),  \alpha (u)) + 2(xy*\alpha (z))*{\alpha}^{2}(u), \nonumber \\
& & \{ y*z, \alpha (x),  \alpha (u)\} = - J_{\alpha}(y*z, \alpha (x),  \alpha (u)) + 2(yz*\alpha (x))*{\alpha}^{2}(u), \nonumber \\
& & \{ z*x, \alpha (y),  \alpha (u)\} = - J_{\alpha}(z*x, \alpha (y),  \alpha (u)) + 2(zx*\alpha (y))*{\alpha}^{2}(u). \nonumber
\end{eqnarray*}
Then
\begin{eqnarray*}
{\circlearrowleft}_{(x,y,z)} \{ x*y, \alpha (z),  \alpha (u)\} &=& - J_{\alpha}(x*y, \alpha (z),  \alpha (u)) - J_{\alpha}(y*z, \alpha (x),  \alpha (u)) \nonumber \\
& & -J_{\alpha}(z*x, \alpha (y),  \alpha (u)) + 2J_{\alpha}(x,y,z)*{\alpha}^{2}(u). \nonumber
\end{eqnarray*}
We know \cite{Issa2} that the identity (4) is equivalent to the identity \\
\par
$ J_{\alpha}(\alpha (x), \alpha (y), x*z) = J_{\alpha}(x,y,z)*{\alpha}^{2}(x)$ \\
\\
(see (2)) defining Hom-Malcev algebras \cite{Yau3}. Then, by (4), 
\par
$ J_{\alpha}(x*y, \alpha (z),  \alpha (u)) = J_{\alpha}(x,z,u)*{\alpha}^{2}(y) + {\alpha}^{2}(x) * J_{\alpha}(y,z,u)$ \\
 \hspace*{4.5truecm}$- 2J_{\alpha}(z*u, \alpha (x),  \alpha (y))$,
\par
$ J_{\alpha}(y*z, \alpha (x),  \alpha (u)) = J_{\alpha}(y,x,u)*{\alpha}^{2}(z) + {\alpha}^{2}(y) * J_{\alpha}(z,x,u)$ \\
 \hspace*{4.5truecm}$- 2J_{\alpha}(x*u, \alpha (y),  \alpha (z))$,
\par
$ J_{\alpha}(z*x, \alpha (y),  \alpha (u)) = J_{\alpha}(z,y,u)*{\alpha}^{2}(x) + {\alpha}^{2}(z) * J_{\alpha}(x,y,u)$ \\
 \hspace*{4.5truecm}$- 2J_{\alpha}(y*u, \alpha (z),  \alpha (x))$. \\
Therefore
\par
$ {\circlearrowleft}_{(x,y,z)} \{ x*y, \alpha (z),  \alpha (u)\} = -2J_{\alpha}(x,z,u)*{\alpha}^{2}(y) - 2{\alpha}^{2}(x) * J_{\alpha}(y,z,u)$ \\ 
\hspace*{5truecm}$-2J_{\alpha}(y,x,u)*{\alpha}^{2}(z) + 2J_{\alpha}(z*u, \alpha (x),  \alpha (y))$ \\ 
\hspace*{5truecm}$ + 2J_{\alpha}(x*u, \alpha (y),  \alpha (z)) + 2J_{\alpha}(y*u, \alpha (z),  \alpha (x))$ \\
\hspace*{5truecm}$ + 2J_{\alpha}(x,y,z)*{\alpha}^{2}(u)$. \\
Now, observe that by the identity \\
\par
$ J_{\alpha}(\alpha (w),  \alpha (y), x*z) + J_{\alpha}(\alpha (x),  \alpha (y), w*z) =$
\par
 $J_{\alpha}(w,y,z)*{\alpha}^{2}(x) + J_{\alpha}(x,y,z)*{\alpha}^{2}(w)$ \hfill (8) \\
\\
(see (8) in \cite{Yau3}) which is shown \cite{Yau3} to be equivalent to (2), we have
\par
$ -2J_{\alpha}(y,x,u)*{\alpha}^{2}(z) + 2J_{\alpha}(y*u, \alpha (z),  \alpha (x)) + 2J_{\alpha}(x,y,z)*{\alpha}^{2}(u)=$
\par
$ 2J_{\alpha}(\alpha (u),  \alpha (x), z*y)$ \\
so that
\par
$ {\circlearrowleft}_{(x,y,z)} \{ x*y, \alpha (z),  \alpha (u)\} = -2J_{\alpha}(x,z,u)*{\alpha}^{2}(y) - 2{\alpha}^{2}(x) * J_{\alpha}(y,z,u) $ \\
\hspace*{5truecm}$ +2J_{\alpha}(\alpha (u),  \alpha (x), z*y) + 2J_{\alpha}(z*u, \alpha (x),  \alpha (y))$ \\
\hspace*{5truecm}$ +2J_{\alpha}(x*u, \alpha (y),  \alpha (z))$. \\
Next, by (4), we have
\par
$ 2J_{\alpha}(z*y, \alpha (u),  \alpha (x)) = -J_{\alpha}(x*u, \alpha (y),  \alpha (z)) + {\alpha}^{2}(x) * J_{\alpha}(y,z,u) $ \\
\hspace*{5.0truecm}$ +{\alpha}^{2}(u) * J_{\alpha}(x,z,y)$,
\par
$ 2J_{\alpha}(z*u, \alpha (x),  \alpha (y)) = -J_{\alpha}(x*y, \alpha (z),  \alpha (u)) + J_{\alpha}(x,z,u)*{\alpha}^{2}(y) $ \\
\hspace*{5.0truecm}$ + {\alpha}^{2}(x) * J_{\alpha}(y,z,u)$ \\
so that
\par
$ {\circlearrowleft}_{(x,y,z)} \{ x*y, \alpha (z),  \alpha (u)\} = -2J_{\alpha}(x,z,u)*{\alpha}^{2}(y) - 2{\alpha}^{2}(x) * J_{\alpha}(y,z,u) $ \\
\hspace*{5.5truecm}$ -J_{\alpha}(x*u, \alpha (y),  \alpha (z)) + {\alpha}^{2}(x) * J_{\alpha}(y,z,u)$ \\
\hspace*{5.5truecm}$ + {\alpha}^{2}(u) * J_{\alpha}(x,z,y) - J_{\alpha}(x*y, \alpha (z),  \alpha (u))$ \\
\hspace*{5.5truecm}$ +J_{\alpha}(x,z,u)*{\alpha}^{2}(y) + {\alpha}^{2}(x) * J_{\alpha}(y,z,u)$ \\
\hspace*{5.5truecm}$ +2J_{\alpha}(x*u, \alpha (y),  \alpha (z)) $ \\
\hspace*{5.5truecm}$ = -J_{\alpha}(x,z,u)*{\alpha}^{2}(y) + {\alpha}^{2}(u) * J_{\alpha}(x,z,y)$ \\
\hspace*{5.5truecm}$ +J_{\alpha}(x*u, \alpha (y),  \alpha (z)) - J_{\alpha}(x*y, \alpha (z),  \alpha (u))$. \\
By (8), we note that
\par
$ J_{\alpha}(x*u, \alpha (y),  \alpha (z)) - J_{\alpha}(x*y, \alpha (z),  \alpha (u))$
\par
$ = -J_{\alpha}(y,z,x)*{\alpha}^{2}(u) -J_{\alpha}(u,z,x)*{\alpha}^{2}(y)$. \\
Therefore, from the last expression of $ {\circlearrowleft}_{(x,y,z)} \{ x*y, \alpha (z),  \alpha (u)\}$ above, we get
\par
$ {\circlearrowleft}_{(x,y,z)} \{ x*y, \alpha (z),  \alpha (u)\} = -J_{\alpha}(x,z,u)*{\alpha}^{2}(y) + {\alpha}^{2}(u) * J_{\alpha}(x,z,y)$ \\
\hspace*{5.5truecm}$-J_{\alpha}(y,z,x)*{\alpha}^{2}(u) -J_{\alpha}(u,z,x)*{\alpha}^{2}(y)$ \\
\hspace*{5.5truecm}$= 0$ \\
and thus {\bf (HLY6)} holds.
\par
The checking of {\bf (HLY7)} is given by the proof of Proposition 2.1 above.
\par
Finally we check the validity of {\bf (HLY8)} for $(A, *, \alpha)$. We have \\
\\
$\{ \{ x,y,u \}, {\alpha}^{2}(v), {\alpha}^{2}(w) \} + \{ {\alpha}^{2}(u), \{ x,y,v \}, {\alpha}^{2}(w) \} + \{ {\alpha}^{2}(u), {\alpha}^{2}(v), \{ x,y,w \} \}$ \\
\par
$= \{ x,y,u \} {\alpha}^{2}(v) * {\alpha}^{3}(w) - {\alpha}^{2}(v) {\alpha}^{2}(w) * {\alpha}(\{ x,y,u \}) - {\alpha}^{2}(w)\{ x,y,u \} * {\alpha}^{3}(v)$
\par
$+ {\alpha}^{2}(u)\{ x,y,v \} * {\alpha}^{3}(w) - \{ x,y,v \}{\alpha}^{2}(w) * {\alpha}^{3}(u) - {\alpha}^{2}(w) {\alpha}^{2}(u) * {\alpha}(\{ x,y,v \})$
\par
$+ {\alpha}^{2}(u){\alpha}^{2}(v) * {\alpha}(\{ x,y,w \}) - {\alpha}^{2}(v)\{ x,y,w \} * {\alpha}^{3}(u) - \{ x,y,w \} {\alpha}^{2}(u) * {\alpha}^{3}(v)$ (by (3)) \\
\par
$=(\{x,y,u\}*{\alpha}^{2}(v) + {\alpha}^{2}(u)*\{x,y,v\})*{\alpha}^{3}(w) + {\alpha}^{2}(u){\alpha}^{2}(v)* \alpha(\{ x, y, w \})$
\par
$-(\{x,y,v\}*{\alpha}^{2}(w) + {\alpha}^{2}(v)*\{x,y,w\})*{\alpha}^{3}(u) - {\alpha}^{2}(v){\alpha}^{2}(w)* \alpha(\{ x,y, u \})$
\par
$-(\{x,y,w\}*{\alpha}^{2}(u) + {\alpha}^{2}(w)*\{x,y,u\})*{\alpha}^{3}(v) - {\alpha}^{2}(w){\alpha}^{2}(u)* \alpha (\{ x,y,v \})$ \\
\par
$=\{ \alpha(x), \alpha(y), u*v \}*{\alpha}^{3}(w)+{\alpha}^{2}(u*v)* \{\alpha (x),\alpha (y),\alpha (w) \}$
\par
$-\{ \alpha(x), \alpha(y), v*w \}*{\alpha}^{3}(u)-{\alpha}^{2}(v*w)* \{\alpha (x),\alpha (y),\alpha (u)\}$
\par
$-\{ \alpha(x), \alpha(y), w*u \}*{\alpha}^{3}(v)-{\alpha}^{2}(w*u)* \{\alpha (x),\alpha (y),\alpha (v)\}$
\par
(by (7) and multiplicativity) \\
\par
$=\{ {\alpha}^{2}(x), {\alpha}^{2}(y), uv*\alpha(w) \}-\{{\alpha}^{2}(x), {\alpha}^{2}(y),vw*\alpha(u)\}$
\par
$- \{ {\alpha}^{2}(x), {\alpha}^{2}(y), wu*\alpha(v) \}$ (by (7)) \\
\par
$=\{ {\alpha}^{2}(x), {\alpha}^{2}(y), \{ u,v,w \} \}$ (by (3)) \\
\\
and thus we get {\bf (HLY8)}.
\par
This completes the proof.\hfill $\square$ \\
\par
The untwisted version of Theorems 1.1 and 1.2 above is proved by K. Yamaguti in \cite{Yam2}, \cite{Yam3}. One observes that Yamaguti's proofs are quite 
different of our proofs above in the situation of Hom-algebras. It could be of some interest to know in which extent the Yamaguti's approach can be applied 
in the Hom-algebra setting.

\end{document}